\documentclass[reqno,12pt]{article}
\usepackage{amsfonts}
\usepackage{graphicx}
\usepackage{epstopdf}
\usepackage{amsmath}
\usepackage{amsthm}
\usepackage{latexsym,bm}
\usepackage{amssymb}
\usepackage{indentfirst}

\usepackage{cite}
\usepackage{setspace}
\newtheorem{thm}{Theorem}[section]

{\theoremstyle{definition}

}

\numberwithin{equation}{section}

\begin{document}
\title{\bf  An organizing principle in the study of the Jacobian Conjecture

}
\author{ {\sc \rm Frederico  Xavier}}
\date{}
\maketitle

\begin{abstract}  \noindent Let $\Omega$ be an irreducible  component of the locus of polynomial  maps 
$  F:\mathbb C^n  \to \mathbb C^n$  satisfying     $\text{\rm deg} F\leq k$ and   $\text{det}DF=1$. It is shown that either   $\Omega \subset \text{\rm Aut} (\mathbb C^n)$, as claimed by the Jacobian Conjecture,    or   the {\em general}   $F\in \Omega$  is   not an automorphism. \end{abstract}

\section{Introduction}

The  Jacobian Conjecture  (JC)  claims that every  polynomial map $F:\mathbb C^n \to  \mathbb C^n$   satisfying    $\text{det}DF=1$     is          invertible \cite{1, E}. Let     $\mathcal P_n^k$ be  the  vector   space  of   all  polynomial maps 
$F=(F_1, \dots, F_n):\mathbb C^n\to \mathbb C^n$  with   $\text{\rm deg} F \stackrel{\text{def}}{=} \max\text{\rm deg}F_j\leq k,$    
and   $V_n^k \subset \mathcal P_n^k$ the algebraic  locus     given    by     $\text{\rm det} DF=1$. 
In this notation JC is the claim that    $\bigcup_{n,k\in \mathbb N}V_n^k\subset \text{Aut}(\mathbb C^n)$.  
Our goal is to   investigate     JC   at the finer level of a fixed     irreducible component  $\Omega$ of $V_n^k$.

 \begin{thm} \label{T1} Either every $F\in \Omega$ is invertible  or   the general $F\in \Omega$ is non-invertible.  
 Thus, either JC holds in   $\Omega$ or \underline{almost every}  $F\in \Omega$ is  actually  a counterexample of JC.  \end{thm}

As    injective polynomial maps  $\mathbb C^n\to \mathbb C^n$  are  bijective (the Ax-Grothendieck theorem), one   obtains  the following  probabilistic criterion  to resolve the ambiguity in the dichotomy: 
{\em \underline{Almost surely}, the correct alternative in Thm.\ref {T1}  is  the first one   if and only if   for a 
\underline {random  choice} of $(F,c)\in \Omega\times \mathbb C^n$ the solution of $F(z)=c$   in $\mathbb C^n$ exists and is unique.} 
    
JC   has been verified  for several   classes of   maps with special properties (e.g. quadratic maps \cite{W}),   but  all  attempts at  establishing  it in full  generality were found to be faulty. The fact that Thm.\ref{T1} is  a dichotomy  puts this history in perspective and opens  up   the possibility that    maybe JC fails for some pair  $(V_n^k, \Omega)$.  It is conceivable that   an application of the  above criterion  to  all  components of all $V_n^k$  may   lead to  counterexamples. 
The silver lining in such a  difficult   task  is that one only needs for the algorithm  to reach  the  conjectural \lq\lq  right\rq\rq  $\Omega$, as   the counterexamples  would then  be  nearly \lq\lq  impossible\rq\rq   to miss.

\section{Proof of Theorem \ref{T1}} \label{sec2}

Because  there is a wider audience for  the Jacobian Conjecture, the conceptual presentation that follows contains more details than it is usually the case  in   research papers.  We begin by fixing notation and terminology. Let  $N(n,k)$ be the dimension of the   complex  vector space of polynomials in $z_1, \dots, z_n$ that have degree at most $k$, so that $ \text{dim}\; \mathcal P_n^k=nN(n,k)$.  A subset    $C$ of  $\mathcal P_n^k$  is   called   algebraic if,  for some (and hence every)   linear  isomorphism $\Lambda: \mathbb C^{nN(n,k)} \to \mathcal P_n^k$,      
$\Lambda^{-1}(C)$ is   algebraic in the  affine sense, i.e. $\Lambda^{-1}(C)$ is the intersection of zero sets of polynomials. 

If we assume   that the nonempty set $I_n^k\stackrel{\text{def}}{=}V_n^k\cap \text{Aut}(\mathbb C^n) $ is  algebraic in  $\mathcal P_n^k$ then the theorem follows   in a  straightforward way. Indeed, if $\Omega \cap I_n^k=\emptyset$   the second alternative in Thm.\ref{T1}  holds in the stronger sense that all maps in $\Omega$ are non-invertible.  Otherwise,   
$\Omega \cap I_n^k\subset\Omega$ 
is     an inclusion of  non-empty algebraic sets.  If    $\Omega \cap I_n^k= \Omega$,    the first alternative  holds. If     
$\emptyset \neq \Omega \cap I_n^k\neq \Omega$,    the irreducibility of $\Omega$ implies   the inequality 
$\text{dim} (\Omega\cap I_n^k)< \text{dim}\Omega$. In this case the second alternative  holds, where the relevant   non-empty Zariski  open (hence dense) subset of 
$\Omega$ where the  general  maps live   is
 $\Omega-(\Omega\cap I_n^k)$.

The proof  that 
  $I_n^k$ is   algebraic in  $\mathcal P_n^k$ will be the object of the remainder of this  note.
To implement  analytic arguments  it is convenient to  work with  the   isomorphism $\Lambda=\Lambda_n^k$  described below.
Let  $\mathcal B$  be the set of   	polynomials  $z_1^{\alpha_1}\dots z_n^{\alpha_n}$,  $\alpha_j\geq 0$, 
$\alpha_1 +\dots + \alpha_n \leq k$, endowed with some order. The  cardinality of $\mathcal B$ is  of course $N(n,k)$.
Define an      isomorphism   $\Lambda_n^k: \mathbb C^{nN(n,k)} \to \mathcal P_n^k$ as follows.  Multiply   each of  the first  $N(n, k)$ coordinates of $c\in \mathbb C^{nN(n,k)}$  by the  respective   polynomial in $\mathcal B$,   and then add  all    terms. The resulting polynomial is, by definition,   the first  component of  $\Lambda_n^kc$. Similarly, using  the second block of $N(n, k)$ coordinates of $c$, re-indexed from $1$ to $N(n,k)$,    multiply each coordinate by the     corresponding   polynomial in $\mathcal B$ and  add all   terms  to   produce   the  second   component of  $\Lambda_n^kc$. Proceed in this fashion,   up to the re-indexed $n$-th block of $N(n,k)$ 
coordinates of $c$, to define all components     of $\Lambda_n^kc$.
For any given  $F\in \mathcal P_n^k$  the vector $(\Lambda_n^k)^{-1}F$  is    simply  an ordered    list of  the coefficients of  all    components $F_i$ of $F$.

The set  $(\Lambda_n^k)^{-1}V_n^k $ is  obviously  Zariski-closed,  since it is    defined by  the polynomial relations among the coefficients of all  $F_i$ that  arise  by setting  $\text{det} DF(z)= 1, \; z\in \mathbb C^n$.   

To prove the theorem we must   argue that 
 $(\Lambda_n^k)^{-1}I_n^k$  is  also Zariski-closed.  In order to do this, consider the set 
 $$ T_j\stackrel{\text{def}}{=} \big\{(c, z, w) \in \mathbb C^{nN(n,k)}\times \mathbb C^n \times \mathbb C^n)  \;  \vert \;  \Lambda_n^k c(z)-\Lambda_n^k c(w)=0, \;    z_j-w_j\neq 0 \big\},$$ 
  where  $1\leq j \leq n$,  $z=(z_1, \dots, z_n),  w=(w_1, \dots, w_n).$  
Relative to the Zariski topology of $\mathbb C^{nN(n,k)+2n} \approx \mathbb C^{nN(n,k)}\times \mathbb C^n \times \mathbb C^n$, $T_j$ is  the intersection of a closed  set and an open  one.      Hence, 
the set $T=T_1\cup\dots\cup T_n$ 
 is  constructible.    
By a   special case of Chevalley's theorem, the   projection $\pi_1(T)\subset \mathbb C^{nN(n,k)}$ is also constructible.  
Manifestly, an  element of  $\pi_1(T)$  is  a    list of  the coefficients of  a  {\em non-injective} map     in $\mathcal P_n^k$.  Conversely, any such map  arises in this way.

Since  the class of constructible sets  is   invariant under  the operations of taking  complements and finite intersections, 
$(\mathbb C^{nN(n,k)}-\pi_1(T))\cap (\Lambda_n^k)^{-1}V_n^k$ 
 is  also  constructible. 
Clearly, an element of    this   set 
 is a list   of   the coefficients of  some    {\em injective}  polynomial map  that belongs   to     $V_n^k$ (and  conversely).
By the Ax-Grothendieck theorem   one   concludes that 
$(\Lambda_n^k)^{-1}I_n^k =(\mathbb C^{nN(n,k)}-\pi_1(T))\cap (\Lambda_n^k)^{-1}V_n^k.$ 
It follows from  this equality and the previous discussion that  $(\Lambda_n^k)^{-1}I_n^k$    is constructible. 
To finish the proof of  Thm.\ref{T1}  we   must   improve matters and 
argue  that  $(\Lambda_n^k)^{-1}I_n^k$
is  actually  Zariski-closed. 
Since     the  Euclidean  and Zariski closures of a constructible set   are  equal,    it suffices to prove  that 
$(\Lambda_n^k)^{-1}I_n^k$   is   classically    closed. 

From  this point on   the arguments are   analytic and topological.
Consider    a  sequence $(c^{(j)})$ in $ \mathbb C^{nN(n,k)}$,  $j\in \mathbb N$,  satisfying the following conditions: 
\vskip5pt
\noindent i) $c^{(j)}\stackrel{j\to \infty}{\longrightarrow} c\in \mathbb C^{nN(n,k)}$. 
\newline 
  ii)   $\text{\rm det} D\Lambda_n^k{c^{(j)}}(z)=1, z\in \mathbb C^n$.   
  \newline
iii) $\Lambda_n^k{c^{(j)}}:\mathbb C^n\to \mathbb C^n$ is  injective. 

\vskip5pt
As observed before, the injective  maps in $\mathcal P_n^k$  are  invertible. Conditions ii) and iii),  taken together,   are then  equivalent to  $c^{(j)}\in (\Lambda_n^k)^{-1}I_n^k$.  Since our immediate goal is to show that $(\Lambda_n^k)^{-1}I_n^k$ is closed in the Euclidean topology, we must argue that  $c\in (\Lambda_n^k)^{-1}I_n^k$ as well, i.e. that ii) and iii)  still   hold  if  $c^{(j)}$ is replaced by $c$.
Regarding ii), continuity and i)  imply     $\text{\rm det} D\Lambda_n^k{c}(z)=1$.  
The analysis of  the analogue of iii) is a bit more involved.
Since   $\text{\rm det} D\Lambda_n^k{c}(z)$ is nowhere zero, the   complex  inverse function theorem   implies that   $\Lambda_n^k c:\mathbb C^n\to \mathbb C^n$  has     local analytic (hence continuous) inverses, and so  its  realification 
$G\stackrel{\text{def}}{=}\widetilde{\Lambda_n^k c}:\mathbb R^{2n}\to \mathbb R^{2n},$
 obtained by taking real and imaginary parts of the components of $\Lambda_n^k{c}$,  is a local homeomorphism. 

Next,   we need a topological observation: 
\vskip5pt
\noindent ($\dagger$) \label{dagger}  If  $(G_j)$ is a sequence of injective local homeomorphisms $G_j: \mathbb R^m \to \mathbb R^m$  that  converges locally uniformly  to a local homeomorphism $G:\mathbb R^m \to \mathbb R^m$, then $G$ is also   injective.
\vskip5pt
\noindent (In passing, elementary one-dimensional examples show that the analogue of ($\dagger$) fails if injectivity is replaced by surjectivity.)  For completeness, we provide    the  arguments for  ($\dagger$).  Assume      by contradiction      that $G(p_1)=G(p_2)=q$,  with $p_1\neq p_2$, and let $U_1, U_2$ be disjoint bounded  neighborhoods of $p_1$ and $p_2$ whose closures  are mapped homeomorphically by $G$ onto a  neighborhood $W$ of $q$.  
As $q\notin G(\partial U_i)$  for   $i=1, 2$,  the topological degree $d(G, U_i, q)$     is  defined and   has value  $\pm 1$ because $\overline{U_i}\stackrel {G}{\to} W$ is a homeomorphism. 
Since  $G_j\vert {\overline {U_i}} \stackrel {j\to \infty} {\longrightarrow} G\vert \overline {U_i}$ uniformly, $d(G_j, U_i, q)$ is also defined for all  sufficiently   large $j$ and 
$d(G_j, U_i, q)=d(G, U_i, q)\neq 0.$
In particular, 
 $G_j(x)=q$ has solutions in  both disjoint sets    $U_1$ and $U_2$,   a contradiction to the injectivity of $G_j$.  
 
 \qed
 
It is easy to see that for any given  compact set $K\subset \mathbb C^n$    there exists  a constant $C(n,k,K)$, independent of $j$,  such that 
 \begin{align} \label{est} & \sup_{z\in K}|| \Lambda_n^kc^{(j)}(z)-\Lambda_n^k c(z)||=\sup_{z\in K} || \Lambda_n^k(c^{(j)}-c)(z)|| 
  \leq C(n,k,K)||c^{(j)}-c||.
  \end{align}

A consequence of  i) and (\ref{est}) is that      $\Lambda_n^kc^{(j)} \stackrel{j\to \infty} {\longrightarrow} \Lambda_n^k c$ locally uniformly.
 In particular, the  induced injective local homeomorphisms 
 $G_j\stackrel{\text{def}}{=}\widetilde{\Lambda_n^k c^{(j)}}:\mathbb R^{2n}\to \mathbb R^{2n}$ 
  converge  locally uniformly  
 to the local homeomorphism $G\stackrel{\text{def}}{=}\widetilde{\Lambda_n^k c}$. 
It follows from ($\dagger$) that $\widetilde{\Lambda_n^k c}$ is injective,  and evidently so is ${\Lambda_n^k c}$.
Thus, the analogue of iii)  with $c$ replacing $c^{(j)}$ holds, and    $(\Lambda_n^k)^{-1}I_n^k$ is classically closed,  completing  the proof of the theorem.

\section{Related results} 
Regarding the    problem in real algebraic geometry of  understanding when   a    polynomial map  $G:\mathbb R^n\to \mathbb R^n$ can be inverted,    in \cite{X2}  we use ideas from nonlinear optimization and the Tarski-Seidenberg theorem to  establish the existence of     a  subset  $\text{\rm ker} (G)$ of $\mathbb R^n$    satisfying:
\vskip5pt
\noindent i) $\text{\rm ker} (G)$ is semialgebraic.   
\vskip1pt
\noindent ii)  $0\in \text{\rm ker} (G)$ and  $\text{\rm ker} (G)=\{0\}$ if and only if $G$ is  invertible.
\vskip1pt
\noindent iii) If $G$ is linear then $\text{\rm ker} (G)$ is the usual kernel.
\vskip5pt
An earlier version of these results, valid  for polynomial local diffeomorphisms, can be found in \cite{X1}.
For  analytic, topological, and geometric mechanisms  of global inversion of maps   in the {\em smooth} category,    we refer to   the semi-expository  paper \cite{X3} and the references therein.

\vskip40pt

\noindent  Frederico Xavier \\
 Professor  Emeritus \\
  University of Notre Dame  \& Texas Christian University \\
Department of Mathematics \\
 fxavier@nd.edu

\end{document}